\documentclass[a4paper,fleqn,oneside]{article}
\usepackage{amsmath}
\usepackage{amsthm}
\usepackage{amsfonts}
\usepackage{amssymb}
\usepackage{subeqnarray}
\usepackage {graphicx}
\usepackage{float}
\usepackage[]{hyperref}
\def\E{{\mathbb E}}
\def\bbE{{\mathbb E}}

\def\S{{\mathbb S}}

\def\bP{{\bf P}}
\def\bE{{\bf E}}

\def\cG{{\mathcal G}}
\def\cF{{\mathcal F}}

\def\cG{{\mathcal G}}

\def\fM{{\mathfrak M}}

\def\vtilde{\tilde v}

\def\wtilde{\tilde w}

\newsavebox{\limn}
\sbox{\limn}{\large $\lim \limits_{n\rightarrow\infty}$}
\newcounter{wspolnynrm}[section]
\newcounter{wspolnynrml}[section]

\newtheorem{Theorem}[wspolnynrm]{{Theorem}}
\newtheorem{Lemma}[wspolnynrml]{Lemma}

\def\rightbox{\protect\vspace*{-2ex}
\begin{flushright}\(\blacksquare\)\end{flushright}}
\newenvironment{Proof}{{\sc Proof.}\hspace{1mm}}{\rightbox}

\begin{document}
\title{Shelf Life of Candidates in the Generalized Secretary Problem.%
\thanks{Final form: 10th of October 2007; File name: \jobname.tex}
}


\author{Krzysztof Szajowski\thanks{Institute of Mathematics, Wroc\l{}w University
of Technology, Wybrze\.ze Wyspia\'nskiego~27, 50-370 Wroc\l{}aw, Poland;
\newline{\it e-mail: szajow@im.pwr.wroc.pl}}
and
Mitsushi Tamaki\thanks{Aichi University, Nagoya Campus: 370 Kurozasa,
Miyoshi, Nishikamo, Aichi 470-0296, Japan; \newline{\it e-mail: tamaki@vega.aichi-u.ac.jp}},
}

  \date{October 10, 2007}

  \maketitle
\begin{abstract}
 A version of the secretary problem called the duration problem, in which the objective is
 to maximize the time of possession of relatively best objects or the second best,
 is treated. It is shown that in this duration problem there are threshold numbers
 $(k_1^\star,k_2^\star)$ such that the optimal strategy immediately selects a
 relatively best object if it appears after time $k_1^\star$ and a relatively second best object
 if it appears after moment $k_2^\star$. When number of objects tends to infinity the thresholds values
 are $\lfloor 0.417188N\rfloor$ and $\rfloor 0.120381N\rfloor$, respectively. The asymptotic mean time of
 shelf life of the object is $0.403827N$.



\end{abstract}

  \section{\label{sec1shell}Introduction and summary}
The duration problem for the classical no-information secretary problem has
been investigated for the first time by Ferguson, Hardwick and Tamaki~\cite{ferhartam92:own}.
It is a sequential selection problem which is a variation of the classical secretary problem (CSP),
treated for example, by Gilbert and Mosteller~\cite{gilmos66}. The aim in CSP is to examine
items ranked from $1$ to $N$ by random selection without replacement, one at a time, and to win is to stop
at any item whose overall (absolute) rank belongs to the given set of ranks (in the basic version this set
contains the rank $1$ only), given only the relative ranks of the items drawn so far.
Since the articles by Gardner~\cite{gar:sci60a,gar:sci60b} the secretary problem has been extended
and generalized in many different directions. Excellent reviews of the
development of this colourful problem and its extensions have been given by
Rose~\cite{rose}, Freeman~\cite{fre}, Samuels~\cite{sam91:secretary} and
Ferguson~\cite{fer89:who}.

The basic form of the duration problem can be described as follows. A set of $N$ rankable objects
appears as in CSP. As each object appears, we must decide
to select or reject it based on the relative ranks of the objects. The
payoff is the length of time we are in possession of a relatively best
object. Thus we will only select a relatively best object, receiving a
payoff of one as we do so and an additional one for each new observation
as long as the selected object remains relatively best.

Though Ferguson, Hardwick and Tamaki~\cite{ferhartam92:own} considered
the various duration models extensively, they confined themselves to the
study of the shelf life the relatively best items. In his seminal paper Gnedin~\cite{gne05:objectives} has shown that such problems are equivalent to the analogous best-choice problems with random horizon N, uniformly
distributed from 1 to n. 

In this paper our goal is different. We attempt to extend the problem to choose and keep the best or the second best items. For simplicity we refer  to a relatively best or a second best object as a candidate.
We receive each time a unit payoff as long as  either of the chosen objects remains a {\bf candidate}.
Obviously only candidates can be chosen, the objective being to
maximize expected  payoff. This problem can be viewed from another
perspective as follows. Let us observe at moment $i$ the relatively second  candidate and let us denote $T(i)$
the time  of the first candidate after time $i$ (\emph{i.e.} the relatively best or the second best item)
if there is one, and $N+1$ if there is none. If we observe at $i$ the relatively best item then $T(i)$
is the moment when new item appears which changes the relative rank of $i$th item to the no candidate rank.
The time $T(i)-i$ is called duration of the candidate selected at time $i$. The objective is to find a stopping
time $\tau^*)$ such that
\begin{equation}\label{problem}
 v_N=E\left[ \frac{T(\tau^*)-\tau^*}{N}\right]
 =\sup\limits_{\tau\in \fM^N} E\left[\frac{T(\tau)-\tau}{N}\right],
\end{equation}
where $\tau$ denotes the stopping time.
\par
This problem will be discussed in Section~\ref{sec2shell}. A Markov chain model will be formulated
and the optimal strategy in Section~\ref{sec22} will be derived. This section is based mainly
on the suggestion from~\cite{dynyush} and the results by~\cite{sza1} and \cite{sucsza02}.
It can be shown that, there exists an optimal threshold stopping time such that, it immediately
selects a best candidate if it appears after or on time $k_1^\star$ and it immediately selects
a second best candidate if it appears after time $k_2^\star>k_1^\star$. In Section~\ref{sec3shell}
we investigate  the asymptotics as $n\rightarrow\infty$. $k_1^\star/n$ proves to converge to
$a\cong 0.120381$ and $k_2^\star/n$ to $b\cong 0.417188$. The asymptotic mean time of shelf life of the relatively
best or the second best object is $0.403827N$.

  \section{\label{sec2shell}Markov model for the shell life of the best and the second best}
The models we consider here are so called no information model where the  decision to select an object
is based only on the relative ranks of the objects observed so far. Let $\S=\{1,2,\ldots,N\}$ be the set
of ranks of items and $\{x_1,x_2,\ldots,x_N\}$ their permutation. We assume that all of them are
equally likely. If $X_k$ is rank of $k$-th candidate we define
\[
Y_k=\#\{1\leq i\leq k: X_i\leq X_k\}.
\]
The random variable $Y_k$ is called {\it relative rank} of $k$-th candidate
with respect of items investigated to the moment $k$.

We observe sequentially the permutation of items from the set $\S$. The
mathematical model of such experiment is the probability space
$(\Omega,\cF,\bP)$. The elementary events are permutations of the elements from
$\S$ and the probability measure $\bP$ is uniform distribution on $\Omega$. The
observation of random variables $Y_k$, $k=1,2,\ldots,N$, generate the sequence
of $\sigma$-fields $\cF_k=\sigma\{Y_1,Y_2,\ldots,Y_k\}$, $k=1,2,\ldots,N$. The
random variables $Y_k$ are independent and $\bP\{Y_k=i\}=\frac{1}{k}$.

Denote by $\fM^N$ the set of all Markov moments $\tau$ with respect to
$\sigma$-fields $\{\cF_k\}_{k=1}^N$. The decision maker observe the stream of relative ranks.
When $Y_i\in A=\{1,2\}$ it is the potential candidate for the absolutly first or second item.
Sometimes it is enough to keep such candidate by some period to get profit which is proportional to
the shell file of candidate. The random variable $T_i$ is defined as the moment when the keeping candidate
stops to be the candidate. It can be described by $Y_s$ for $s=i,i+1,\ldots,T_i$.

\subsection{\label{TiDistribution}Distribution of $T_i$}
There are two cases:
\begin{description}
\item{$\mathbf{Y_i=2:}$} in this case $T_i=k$ when $Y_i=2,Y_{i+1}>2,Y_{i+2}>2,\ldots,Y_{k-1}>2,Y_k\in A$.
We have for $i<k\leq N$
\begin{subeqnarray}\label{rel2timed}
\slabel{rel2timedk}
\bP\{T_i=k|Y_i=2\}&=&\bP\{Y_i=2,Y_{i+1}>2,Y_{i+2}>2,\ldots,Y_{k-1}>2,Y_k\in A|Y_i=2\}\\
      \nonumber&=&\frac{i-1}{i+1}\frac{i}{i+2}\ldots\frac{k-3}{k-1}\frac{2}{k}=\frac{2(i-1)i}{(k-2)(k-1)k};
\text{and}\\
\slabel{rel2timedN}
\bP\{T_i=N+1|Y_i=2\}&=&\bP\{Y_i=2,Y_{i+1}>2,Y_{i+2}>2,\ldots,Y_{N-1}>2,Y_N>2|Y_i=2\}\\
\nonumber&=&1-\sum_{s=i+1}^N\frac{2(i-1)i}{(s-2)(s-1)s} =\frac{i(i-1)}{N(N-1)}.
\end{subeqnarray}

\item{$\mathbf{Y_i=1:}$} the random variable $T_i=k$ if there exists $s\in\{i+1,\ldots,k-1\}$ such that
 $Y_i=1,Y_{i+1}>1,Y_{i+2}>1,\ldots,Y_{s-1}>1,Y_s=1,Y_{s+1}>2,\ldots,Y_{k-1}>2,Y_k\in A$. We have for $i<k\leq N$
\begin{subeqnarray}\label{rel1timed}
\slabel{rel1timedk}
\bP\{T_i=k|Y_i=1\}&=&
\bP\left\{\bigcup_{s=i+1}^{k-1}\{Y_i=1,Y_{i+1}>1,Y_{i+2}>1,\ldots,Y_{s-1}>1,Y_s=1,\right.\\
                   \nonumber&&\mbox{}\left. Y_{s+1}>2,\ldots,Y_{k-1}>2,Y_k\in A\}|Y_i=1\right\}\\
\nonumber&=&\sum_{s=i+1}^{k-1}\frac{i}{(s-1)s}\frac{2(s-1)s}{(k-2)(k-1)k}=\frac{2i(k-i-1)}{(k-2)(k-1)k}\\
\vspace{-3em}\nonumber\text{and}&&\\
\slabel{rel1timedN}
\bP\{T_i=N+1|Y_i=1\}&=&1-\bP\left\{\bigcup_{s=i+1}^{N}\{Y_i=1,Y_{i+1}>1,\ldots,Y_{s-1}>1,Y_s=1,\right.\\
\nonumber&&\mbox{}\left.Y_{s+1}>2,\ldots,Y_{k-1}>2,Y_N>2\}|Y_i=1\right\}\\
\nonumber&=&1-\sum_{k=i+1}^{N}\sum_{s=i+1}^{k-1}\frac{i}{(s-1)s}\frac{2(s-1)s}{(k-2)(k-1)k}\\
\nonumber&=&1-\sum_{k=i+1}^{N}\frac{2i(k-i-1)}{(k-2)(k-1)k}=1-\frac{(N-i)(N-i-1)}{N(N-1)}=\frac{2Ni-i^2-i}{N(N-1)}
\end{subeqnarray}
\end{description}

The solution of the problem (\ref{problem}) will be performed by its change to the optimal stopping problem
for the embedded Markov chain.
\subsection{\label{sec22}The optimal stopping problem for the embedded Markov chain}
Let us observe that for any $\tau\in\fM^N$
\begin{eqnarray*}\label{problem1}
 \bE\left[\frac{T(\tau)-\tau}{N}\right]
 &=& \sum_{r=1}^N\int_{\{\tau=r\}}\frac{T_r-r}{N}d\bP=
     \sum_{r=1}^N\int_{\{\tau=r\}}\bE\{\frac{T_r-r}{N}|\cF_r\}d\bP\\
 &=& \sum_{r=1}^N\int_{\{\tau=r\}}\bE\{\frac{T_r-r}{N}|Y_r\}d\bP
  = \bE \varphi(\tau,Y_\tau).
\end{eqnarray*}
In the following lemma the function $\varphi(\cdot)$ is calculated. The final form of it is using the
the digamma function ($\digamma$-function) $\psi_n(z)$ (see Abramowitz and Stegun \cite{abrste72:functions} p. 260).
For $n=0$ we will use denotation $\psi(z)$.
This function is defined as $n$th logarithmic derivative of the Euler gamma function $\Gamma(z)$
\begin{eqnarray*}
\psi_n(z)&=&\frac{d^{n+1}}{d\!z^{n+1}}\ln\Gamma(z)\\
         &=&\frac{d^n}{d\!z^n}\frac{\Gamma^{'}(z)}{\Gamma(z)}.
\end{eqnarray*}

\begin{Lemma}\label{thepayoff}
The payoff function $\varphi(k,r)$ has the form
{\renewcommand{\arraystretch}{2}
\begin{equation}\label{varphi}
\varphi(k,r)=\left\{\begin{array}{ll}
\frac{k}{N^2}\left(1+k-N-2N\psi(k)+2N\psi(N)\right)&\mbox{ for $r=1$,}\\
\frac{k(N-k+1)}{N^2}&\mbox{ for $r=2$,}\\
0&\mbox{ otherwise.}
\end{array}
\right.
\end{equation}
}
\end{Lemma}
\begin{Proof}
Based on the distribution of the random variable $T_k$ and the equality $\psi(p+1)-\psi(p)=\frac{1}{p}$ for
the digamma function  we get

\begin{subeqnarray}
\label{varphi1}
\slabel{varphi1A}
\varphi(k,1)&=&\bE\{\frac{T_k-k}{N}|Y_k=1\}=\frac{1}{N}\left(\sum_{s=k+1}^{N+1}(s-k)\bP\{T_i=s|Y_k=1\}\right)\\
\nonumber&=&\frac{1}{N}\left(\sum_{s=k+1}^N\frac{2k(s-k-1)(s-k)}{s(s-1)(s-2)}+(N+1-k)\frac{2Nk-k^2-k}{N(N-1)}\right)\\
\nonumber&=&\frac{k}{N^2}\left(1+k-N-2N(\psi(k)-\psi(N))\right)\\
\slabel{varphi1B}
\varphi(k,2)&=&\bE\{\frac{T_k-k}{N}|Y_k=2\}=\frac{1}{N}\left(\sum_{s=k+1}^{N+1}(s-k)\bP\{T_k=s|Y_k=2\}\right)\\
\nonumber&=&\frac{1}{N}\left(\sum_{s=k+1}^N\frac{2k(k-1)(s-k)}{s(s-1)(s-2)}+(N+1-k)\frac{k(k-1)}{N(N-1)}\right)
=\frac{k(N-k+1)}{N^2}
\end{subeqnarray}
\end{Proof}

\subsubsection{Recursive algorithm}
Let $\fM^N_r=\{\tau\in\fM^N: r\leq \tau\leq N\}$ and
$\wtilde_N(r)=\sup_{\tau\in\fM^N_r}\bE \varphi(\tau,Y_\tau)$. The following algorithm
allows to construct the value of the problem $v_N=w_N(1,1)$.
\begin{equation}\label{step1}
  \wtilde_N(N)=\bE \varphi(N,Y_N)=\frac{2}{N}.
\end{equation}
Let
\begin{subeqnarray}\label{steps}
\slabel{step2}
  w_N(N,r)&=&=\left\{
            \begin{array}{ll}
                  1, &\mbox{ if $r\in A$},\\
                  0, &\mbox{ otherwise,}
            \end{array}
                  \right.\\
\slabel{step3}%
  w_N(k,r)&=&\max\{\varphi(k,r),\bE w_N(r+1,Y_{r+1})\},\\
\slabel{step4}%
  \wtilde_N(k)&=&\bE w_N(k,Y_k)=\frac{1}{k}\sum_{r=1}^k w_N(k,r).
\end{subeqnarray}
We have then $v_N=\wtilde_N(1)$. The optimal stopping time $\tau^*$ is defined
as follows: one have to stop at the first moment $k$ when $Y_k=r$, unless
$w_N(k,r)>\varphi(k,r)$. We can define the stopping set $\Gamma=\{(k,r): \varphi(k,r)\geq \wtilde_N(k)\}$.

\subsubsection{Embedded Markov chain}
Let $a=\max(A)$. The function $\varphi(k,r)$ defined in (\ref{varphi1}) is equal to
$0$ for $r>a$ and non-negative for $r\leq a$. It means that it is rational to choose item for keeping
at moment $k$ when the state $(k,r)$ such that $r\leq a$.
Define $W_0=(1,Y_1)=(1,1)$, $\gamma_t=\inf\{r>\gamma_{t-1}:Y_r\leq
\min(a,r)\}$ ($\inf\emptyset=\infty$) and $W_t=(\gamma_t,Y_{\gamma_t})$.
If $\gamma_t=\infty$ then define $W_t=(\infty,\infty)$. $W_t$ is  the
Markov chain with the state space $\bbE=\{(s,r):s\in\{1,2,\ldots,N\},r\in A\}\cup\{(\infty,\infty)\}$
and the following one step transition probabilities
(see~\cite{sza1})
\begin{equation}\label{trprob}
  \begin{split}
  p(r,s)&=\bP\{W_{t+1}=(s,l_s)|W_t=(r,l_r)\} 
     =\begin{cases}
      \frac{1}{s},             &\text{ if $r<a$, $s=r+1$},\\
      \frac{(r)_a}{(s)_{a+1}}, &\text{ if $a\leq r< s$},\\
      0,                       &\text{ if $r\geq s$ or $r<a$, $s\neq r+1$,}
    \end{cases}
  \end{split}
\end{equation}
with  $p(\infty,\infty)=1$, $p(r,\infty)=1-a\sum_{s=r+1}^Np(r,s)$, where
$(s)_a=s(s-1)(s-2)\ldots(s-a+1)$, $(s)_0=1$. We denote $T\varphi(k,r)=\bE_{(k,r)} \varphi(W_1)$
the mean operator for the function $g:\bbE\rightarrow \Re$. Let
$\cG_t=\sigma\{W_1,W_2,\ldots,W_t\}$ and $\tilde{\fM}^N$ be the set of stopping
times with respect to $\{\cG_t\}_{t=1}^N$. Since $\gamma_t$ is increasing, then
we can define $\tilde{\fM}^N_{r+1}=\{\sigma \in \tilde{\fM}^N : \gamma_\sigma>r\}$.

Let $\bP_{(k,r)}(\cdot)$ be probability measure related to the Markov chain
$W_t$, with trajectory starting in state $(k,r)$ and $\bE_{(k,r)}(\cdot)$ the
expected value with respect to $\bP_{(k,r)}(\cdot)$. From (\ref{trprob}) we can
see that the transition probabilities do not depend on relative ranks, but only
on moments $k$ where items with relative rank $r\leq \min(a,k)$ appear. Based
on the following lemma we can solve the problem~(\ref{problem}) with gain
function (\ref{varphi}) using the embedded Markov chain $\{W_t\}$.
\begin{Lemma} (see \cite{sza1})
  \begin{equation}\label{basic1}
  \bE w_N(k+1,Y_{k+1})=\bE_{(k,r)} w_N(W_1) \mbox{ for every $r\leq \min(a,k)$.}
  \end{equation}
\end{Lemma}

\subsubsection{Solution of the optimal shelf life problem}
First of all the form of $T\varphi(k,r)$ for $(k,r)\in \bbE$ will be given.
\begin{Lemma}\label{exppayoff}
The expected payoff of the fuction $\varphi(\cdot)$ with relation to the embedded Markov chain
$(W_t,\cG_t,\bP_{(1,1)})_{t=0}^N$ has the folloing form:
  \begin{equation}\label{basic2}
  T\varphi(k,r)=\frac{(N-k)((2N-1)k+N-1)}{N^2(N-1)}+2\frac{k}{N^2}(\psi(N)-\psi(k)).
  \end{equation}
\end{Lemma}
\begin{Proof}
The definition of the embedded Markov chain (\ref{trprob}) and the payoff function $\varphi(\cdot)$ in
the lemma \ref{thepayoff} give
\begin{eqnarray*}
T\varphi(k,r)&=&\sum_{j=k+1}^N\sum_{r=1}^2p(k,j)\varphi(j,r)\\
             &=&\sum_{j=k+1}^N\frac{k(k-1)}{j(j-1)(j-2)}\left(\frac{j(2N(\psi(N)-\psi(j))+N-j-1)}{N^2}+
             \frac{j(N-j+1)}{N^2}  \right)  \\
             &=&\sum_{j=k+1}^N\frac{k(k-1)}{j(j-1)(j-2)}
             \left(\frac{j(N-j-1)}{N^2}(j+N+\frac{2N}{j}-\frac{2N}{N-1})  \right)\\
             &=&\frac{(N-k)((2N-1)k+N-1)}{N^2(N-1)}+2\frac{k}{N^2}(\psi(N)-\psi(k)).
\end{eqnarray*}
\end{Proof}

\begin{table}[htb]
\caption{\label{texact1}Decision points and values of the problem}
\begin{center}
\begin{tabular}{||r||c|c|c||}\hline
$N$ &${s_N^1}^\star$&${s_N^2}^\star$&$v_N$\\ \hline
10  &1       &4       &0.527526    \\
20  &2       &8       &0.464357    \\
30  &3       &12      &0.442977    \\
40  &4       &16      &0.432325    \\
50  &6       &21      &0.426411    \\
60  &7       &25      &0.422846    \\
70  &8       &29      &0.420142    \\
80  &9       &33      &0.418024    \\
90  &10      &37      &0.416322    \\
100 &12      &41      &0.415064    \\
200 &24      &83      &0.409431    \\
500 &60      &208     &0.406064    \\
1000&120     &417     &0.404944    \\
\cline{1-4}
$\infty$&[0.120381N] &[0.417188N]&0.403827   \\
\hline\hline
\end{tabular}
\end{center}
\end{table}

Let us denote $A_k(r)=\{(s,r):s>k\}$.
\begin{Theorem}\label{expectedpayoff}
There are constants $k_1^\star$ and $k_2^\star$ such that the optimal stopping time for the problem (\ref{problem})
has the form
\begin{equation*}
      \tau^*=\inf\{t:W_t\in A_{k_1^\star}\cup A_{k_2^\star}\}.
\end{equation*}
The value function
\begin{eqnarray*}
\vtilde_N(k_1^\star,k_2^\star)&=&
\frac{(N (3N-4)-3)+k_1^\star(N-3)\psi(k_1^\star)+k_1^\star\left(2 (N^2-1)\left(\psi_{1}(k_2^\star+1)%
            -\psi_{1}(k_1^\star + 1)\right)\right)}{(N-1) N}\\
&&\mbox{}+\frac{k_1^\star \left(2 (N-1) \psi(N)+(5-3 N) \psi(k_2^\star)\right)}{(N-1) N }\\
&&\mbox{}-\frac{k_1^\star\left(3 N^3+(2 k_2^\star-3) N^2-2 %
          \left({k_2^\star}^2+k_2^\star+2\right)N+{k_2^\star}^2+k_2^\star\right)}%
   {(N-1) N^2 k_2^\star}
\end{eqnarray*}
\end{Theorem}
\begin{Proof} The payoff function $\varphi(\cdot,r)$ for $r\in A$ are unimodal. It can be seen by analysis
of the differences $\varphi(k+1,1)-\varphi(k,1)$ which is decreasing for $k\leq N-1$. The compare of events
related to $T_k=j$ on $Y_k=1$ and $Y_k=2$ leads to the conclusion that $\varphi(k,1)\geq \varphi(k,2)$
for $k\in\{1,2,\ldots,N\}$. The value function $\wtilde_(k)$ is noincreasing by the fact of decrasing number of
stopping times in $\fM^N_k$. At $k=N-1$ both payoff functions are greater than $\wtilde_N(N-1)$.
Let us be $k_2^\star=\inf\{1\leq k\leq N: T\varphi(k,i)\leq \varphi(k,2)\}-1$. We have for $k>k_2^\star$ and $r=1,2$
that $w_N(k,r)=\varphi(k,r)$ and $\wtilde_N(k)=T\varphi(k,r)$.
Let us denote $k_1^\star=\inf\{1\leq k\leq k_2^\star: \wtilde_N(r)<\varphi(k,1) \}$, where
$\wtilde_N(k)=\vtilde(k,k_2^\star)$ and for $k<s$ we have
\begin{eqnarray*}
\vtilde_N(k,s)&=&\sum_{j=k+1}^s\frac{k}{j(j-1)}\varphi(j,1)+\frac{k}{s}\wtilde_N(s)\\
  &=&\frac{(N (3 N-4)-3)+k (N-3) \psi(k) +  k \left(2 (N-1) \psi(N)+(5-3 N) \psi(s)\right)} {(N-1) N}\\
   &&\mbox{}+\frac{ k\left(2 (N^2-1)\left(\psi_1(s+1)-\psi_1(k+1)\right)\right) }{(N-1) N}\\
   &&\mbox{}-\frac{k \left(3 N^3+(2 s-3) N^2-2 \left(s^2+s+2\right) N+s^2+s\right)}%
   {(N-1) N^2 s}.
\end{eqnarray*}
\end{Proof}

The numerical examples of the solution of the shelf life problem with different horizon are given in the table
\ref{texact1}.

  \section{\label{sec3shell}Asymptotic duration problem}

\subsection{\label{sec31}The limes of the finite horizon problem}
Let the number of candidates goes to infinity. For such large number of
candidates we can find optimal solution of (\ref{problem}) based on the following consideration.
As $N\rightarrow\infty$ such that $\frac{r}{N}\rightarrow x\in(0,1]$ the
embedded Markov chain $(W_t,\cF_t,\bP_{(1,1)})$ with state space
$\E=\{1,2,\ldots,N\}\times\{1,2,\ldots,\max(A)\}$ can be treated as Markov
chain $(W^{'}_t,\cF_t,\bP_{(\frac{1}{N},1)})$ on
$\{\frac{1}{N},\frac{2}{N},\ldots,1\}\times\{1,2,\ldots,\max(A)\}$.
\begin{Lemma}\label{limvarphi}
The gain function $\varphi([Nx],r)$ has limit
\begin{equation*}\label{asymgain}
  \varphi(x,r)= \left\{\begin{array}{ll}
     x^2-2 \log (x) x-x&\mbox{ for $r=1$,}\\
     x(1-x)&\mbox{ for $r=2$.}
     \end{array}
     \right.
\end{equation*}
\end{Lemma}
\begin{Proof}
Let us limit the formula $\varphi(k,r)$ given in (\ref{varphi1A}). We get
\begin{subeqnarray}\label{limvarphi2}
\slabel{limvarphi2A}
\varphi(x,1)&=&\lim_{\stackrel{\frac{k}{N}\rightarrow x}{N\rightarrow\infty}}\varphi(k,1)\\
\nonumber&=&\lim_{\stackrel{\frac{k}{N}\rightarrow x}{N\rightarrow\infty}}
\frac{1}{N}\left(\sum_{s=k+1}^N\frac{2k(s-k-1)(s-k)}{s(s-1)(s-2)}+(N+1-k)\frac{2Nk-k^2-k}{N(N-1)}\right)\\
\nonumber&=&\int_{x}^1\frac{2x(z-x)(z-x)}{z^3}d\!z+(1-x)(2x-x^2)\\
\nonumber&=&-2\log(x)-x+x^2\\
\text{and}
\slabel{limvarphi2B}
\varphi(x,2)&=&\lim_{\stackrel{\frac{k}{N}\rightarrow x}{N\rightarrow\infty}}\varphi(k,2)\\
\nonumber&=&\lim_{\stackrel{\frac{k}{N}\rightarrow x}{N\rightarrow\infty}}
\frac{k(N-k+1)}{N^2}=x(1-x).
\end{subeqnarray}
\end{Proof}

We get also
\begin{math}
  \lim_{N\rightarrow\infty}\bE_{(\frac{k}{N},r)}\varphi(W_1)=\bE_{(x,r)}\varphi(W^{''}_1),
\end{math}
where $(W^{''}_t,\cF_t,\bP_{(x,1)})$ is Markov chain with the state space $(0,1]\times\{1,2,\ldots,\max(A)\}$
and the transition density function
\begin{equation}\label{asymprob}
  p(x,y)=\left\{
  \begin{array}{ll}
    \frac{x^a}{y^{a+1}},&\mbox{ $0<x<y\leq 1$,}\\
    0, &\mbox{ $x\geq y$.}
  \end{array}
  \right.
\end{equation}
The expected value with respect to the conditional distribution given in (\ref{asymprob}) is following
\begin{equation}\label{asymexp}
  \bE_{(x,r{'})}\varphi (W_1^{''})=\sum_{r=1}^{\max(A)}\int_x^1p(x,y)\varphi(y,r)dy.
\end{equation}
The recursive formulae (\ref{step1})--(\ref{step4}) in asymptotic case have the
form
\begin{subeqnarray}\label{asymstep}
\slabel{asymstep1}
  v(1)&=&0\\
  \slabel{asymstep2}
  w(x,r)&=&\max\{\varphi(x,r),\bE_{(x,r)}w(W_1^{''})\},\\
  \slabel{asymstep3}
  v(x)&=&\bE_{(x,r)}w(W_1^{''}).
\end{subeqnarray}
The function $w(x,r)$ is the limit of $w_N(k,r)$, when $\frac{k}{N}\rightarrow x\in
(0,1]$, {\it i.e.} $\lim_{N\rightarrow\infty}w_N([Nx],r)=w(x,r)$. The
asymptotic solution we get by `recursive' method based on
(\ref{asymstep1})--(\ref{asymstep3}).

\subsection{\label{sec32}Asymptotic solution of shelf life problem}
 It is of interest to investigate the asymptotic behaviors of
 $k_1^\star$, $k_2^\star$ and $v_N$ as $N$ tends to infinity. The algorithm presented in the section \ref{sec31}
 is used. Based on the lemma \ref{exppayoff} we get
\begin{Lemma}\label{constb}
$\lim_{\stackrel{\frac{k}{N}\rightarrow x}{N\rightarrow\infty}} T\varphi(k,r)=2*(x^2-x-x\log(x))$ and
\begin{equation}\label{constbeq}
b=\lim_{N\rightarrow\infty} \frac{k_2^\star}{N}=-\frac{2}{3}W(-\frac{3}{2}\exp(-\frac{3}{2}))\cong 0.417188,
\end{equation}
where $W(\cdot)$ is the Lambert $W$--function\footnote{This function was introduce by Euler \cite{eul1783:lambert}
with relation to the Lambert transcedental function investigated by Lambert \cite{lam1758:variae}}
(cf. P\'olya and Szeg\H{o} \cite{polsze25:problems})
\end{Lemma}

\begin{Proof}
From the lemma \ref{exppayoff} we get easily the limit of $T\varphi(k,r)$. It allows to formulate the equation
which the $b=\lim_{N\rightarrow\infty} \frac{k_2^\star}{N}$ should fulfil. It is
\begin{eqnarray*}
-2(x\log(x)+x-x^2)=x(1-x).
\end{eqnarray*}
After simple algebra we get that $b$ should fulfil the equation
$-\frac{3}{2}\exp(-\frac{3}{2})=-\frac{3}{2}b\exp(-\frac{3}{2}b)$. The inverse function to the function $h(x)=we^w$
is the Lambert $W$-function. It gives the solution (\ref{constbeq}).
\end{Proof}

\begin{figure}[htbp]
\centerline{\includegraphics[width=3.54in,height=2.18in]{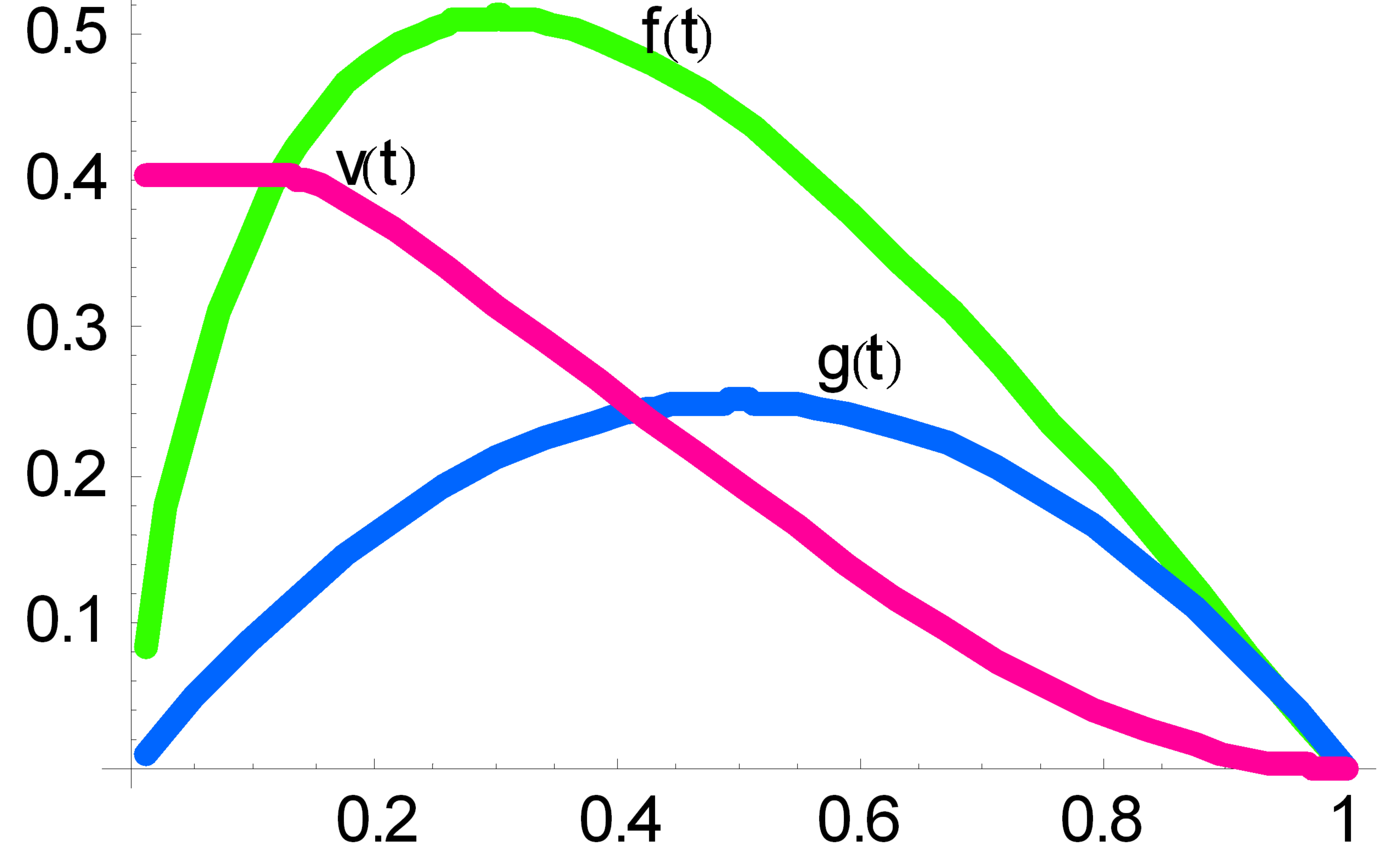}}
\label{fig1}
\end{figure}
  
  \section{\label{sec4shell}Final remarks}
 Thus far we have implicitly assumed that the object, once
 chosen, are possessed until the process terminates. It is possible to extend
 solution for multiple choice duration problem, silmilarly as in \cite{tampeasza98:duration}.
 It will be subject of further research. It is also unknown to the authors if the solution of this
 duration problem has similar solution as the related best choice problem with rundom number of
 objects available. Such coincidence has place for the best candidate duration problem and the
 best choice problem with the random number of objects avalilable considered by Presman and Sonin \cite{preson72}.
 This observation has been given in \cite{fehata92:dur}.



\end{document}